%% file: aaai23.tex
\pdfoutput=1

\documentclass[letterpaper]{article} 
\usepackage{aaai23}  
\usepackage{times}  
\usepackage{helvet}  
\usepackage{courier}  
\usepackage[hyphens]{url}  
\usepackage{graphicx} 
\usepackage{natbib}  
\usepackage{caption} 
\usepackage{algorithm}
\usepackage{algorithmic}

\usepackage{newfloat}
\usepackage{listings}
\DeclareCaptionStyle{ruled}{labelfont=normalfont,labelsep=colon,strut=off} 
\floatstyle{ruled}
\newfloat{listing}{tb}{lst}{}
\floatname{listing}{Listing}
\title{ \textsc{LiMIP} : Lifelong Learning to Solve Mixed Integer Programs}

\author {
    Sahil Manchanda and 
    Sayan Ranu
}
\affiliations {
    Department of Computer Science and Engineering \\ 
    Indian Institute of Technology, Delhi\\
    \{sahil.manchanda,sayanranu\}@cse.iitd.ac.in
}

\usepackage{bibentry}
\input{preamble}

\begin{document}
\maketitle
\begin{abstract}
\input{01_abstract}
\end{abstract}

\input{02_Intro}

\input{04_Method}

\input{05_Experiments}

\input{06_Conclusion}
\input{07_ack}

\bibliography{aaai23}

\appendix
\input{Appendix}

\end{document}

%% file: preamble.tex
\usepackage{soul}
 \usepackage{epstopdf}
\usepackage[utf8]{inputenc}
\usepackage{url,textcomp}  
\usepackage{tablefootnote}
\usepackage{mdwlist}
\usepackage{tabularx}
\usepackage{amsmath,amsfonts,amsthm,bm}
\usepackage{graphicx}
\usepackage{enumitem}
\usepackage{array}
\usepackage{arydshln}
\usepackage{mathtools}
\usepackage{xspace}
\usepackage{xcolor}
\usepackage{graphicx}
\usepackage{amsfonts}
\usepackage{booktabs}

\usepackage{stmaryrd}

\usepackage{todonotes}
\usepackage[toc,page]{appendix}
\usepackage{diagbox}
\usepackage{multirow}
\usepackage{array}
\usepackage{makecell}
\usepackage{fmtcount} 
\usepackage{array,multirow}
\usepackage{tablefootnote}
\usepackage{makecell}
\usepackage{algorithm}
\usepackage{algorithmic}

  \usepackage{graphicx}

\usepackage{subfig}
\usepackage{caption}

 \setlist{nolistsep,leftmargin=*}

\usepackage{amsfonts}

 \newtheorem{defn}{\textbf{Definition}}
  \newtheorem{prob}{\textbf{Problem}}

\def\namemodel{\textsc{LiMIP}\xspace}

\def\GCNN{\textsc{GCNN}\xspace}

%% file: 01_abstract.tex
Mixed Integer programs (MIPs) are typically solved by the Branch-and-Bound algorithm. Recently, \textit{Learning to imitate} fast approximations of the expert \textit{strong branching} heuristic has gained attention due to its success in reducing the running time for solving MIPs. However, existing learning-to-branch methods assume that the entire training data is available in a single session of training. This assumption is often not true, and if the training data is supplied in continual fashion over time, existing techniques suffer from \textit{catastrophic forgetting}. In this work, we study the hitherto unexplored paradigm of \textit{Lifelong Learning to Branch} on Mixed Integer Programs. 
To mitigate catastrophic forgetting, we propose \namemodel, which is powered by the idea of modeling an MIP instance in the form of a bipartite graph, which we map to an embedding space using a bipartite Graph Attention Network. This rich embedding space avoids catastrophic forgetting through the application of knowledge distillation and elastic weight consolidation, wherein we learn the parameters key towards retaining efficacy and are therefore protected from significant drift. We evaluate \namemodel on a series of NP-hard problems and establish that in comparison to existing baselines, \namemodel is up to $50\%$ better when confronted with lifelong learning.

%% file: 02_Intro.tex
\section{Introduction and Related Work}
\looseness=-1
\textit{Combinatorial optimization (CO)} is a subclass of optimization problems that deals with optimizing a certain objective function by selecting a subset of elements from a finite set. Although CO problems are generally NP-hard~\cite{taha2014integer} from a complexity theory viewpoint, still they are solved routinely in diverse fields such as capacity planning, resource allocation, scheduling, and manufacturing~\cite{taha2014integer,gcomb,manchanda2022generalization, benchmarking,greed,neuromlr}. It is common  to formulate most of these CO problems as Mixed integer programs (MIPs)~\cite{achterberg2007constraint}. However, these MIPs are difficult to solve due to the non-convexity of their feasible region. Instead of solving it directly, their LP-relaxed versions are solved~\cite{achterberg2007constraint}. Modern solvers such as SCIP~\cite{gamrath2020scip} employ the Branch-and-Bound (B\&B)~\cite{achterberg2007constraint} algorithm to solve these MIPs. B\&B recursively partitions the solution space into a search tree and then prunes subtrees that provably cannot generate an optimal solution. This is an iterative process, which consists of making sequential decisions such as node selection, branching variable selection, etc., to direct the search procedure. The efficiency of the B\&B algorithm mainly depends upon the branching variable selection and node selection~\cite{achterberg2007constraint}. In this work, we focus on the former. In modern MIP solvers~\cite{gamrath2020scip}, the branching variable decisions are generally based upon hard-coded heuristics designed by experts to direct the search process to solve an MIP~\cite{achterberg2007constraint}. Among various heuristics, \textit{strong branching} is one such heuristic that is highly effective in reducing the size of Branch-and-Bound tree~\cite{achterberg2007constraint, gasse2019exact}. However, its main disadvantage is the extremely high computational cost associated with choosing the best variable to branch. Hence, it is rarely used in practice~\cite{achterberg2007constraint,huang2021branch}. 

\begin{figure}[t!]
\centering
\includegraphics[scale=0.35]{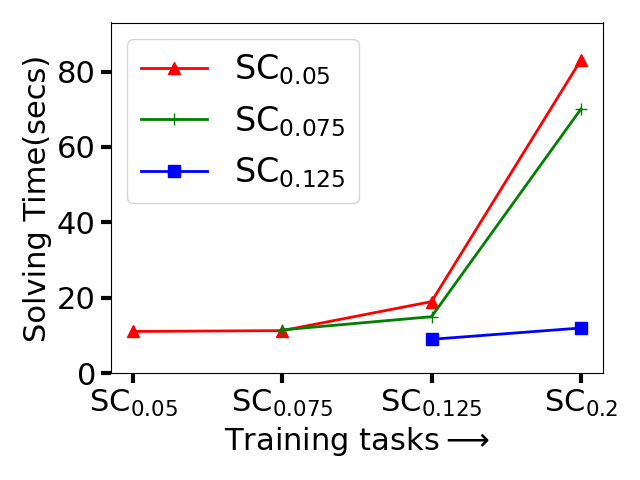}
\caption{Illustration of catastrophic forgetting on SoA method~\cite{gasse2019exact} on  Set Cover dataset. The plot shows the rapid increase in the  geometric mean of solving time on instances of past problems   when the model is further trained \textit{sequentially} on instances of a new problems.  }
\label{fig:intro_forget_sc}
\end{figure}

In order to utilize the advantage of the powerful strong branching heuristic, however at a lower computational cost, recently, multiple algorithms have been developed~\cite{gasse2019exact,khalil2016learning, nair2020solving, gupta2020hybrid}. At their core, they use \textit{imitation learning} to learn \textit{fast approximations} of the \textit{strong branching} heuristic on a family of MIP instances. These algorithms estimate scores of branching candidate variables quickly on unseen but similar MIP instances. This paradigm is popularly known as \textit{learning to branch}. These methods have obtained significant gains in terms of problem solving time over modern MIP solvers such as SCIP.

Despite significant success, existing techniques on learning to branch are limited by the assumption that the entire training data is available in single session of training. 

This assumption is not realistic in the context of MIPs as the semantics of CO problems may keep changing with time. Consequently, the data to train the \textit{learning-to-branch} model is dynamic with updates arriving  sequentially over time. The above scenario is commonly observed in  industries such as shipping and food delivery where entities, such as service locations and warehouses, get added/removed over time. Further, the semantics of the problem such as the demand, supply, customer distribution, facility constraints etc. also fluctuate. Hence, to tackle such scenarios, the model should be capable of learning in an incremental fashion as more data appears over time as retraining the entire model again can be computationally expensive.

The ability to continually learn over time is referred to as \textit{lifelong/continual learning}~\cite{parisi2019continual}. This aligns with the ability of humans to continually acquire skills throughout their lifespan. Lifelong learning aspires to gain more knowledge sequentially and improve existing model as more data arrives.  However, the quintessential failure model of lifelong learning on neural models is \textit{catastrophic forgetting} i.e when new concepts are learned sequentially, the neural model forgets the concepts it learned previously~\cite{kirkpatrick2017overcoming,parisi2019continual}\footnote{Detailed related work present in Appendix A.1 at \url{https://arxiv.org/abs/2208.12226}}.

\looseness=-1
In this context, we analyze the performance of the state-of-the-art learning-to-branch technique \GCNN~\cite{gasse2019exact} in the lifelong learning scenario, i.e., training the model sequentially on different problems. In Fig.~\ref{fig:intro_forget_sc} we observe that the solving time on instances of the  SetCover (\textsc{\textit{SC}}) problem with edge-probability $0.05$, i.e., $\textsc{\textit{SC}}_{0.05}$ increases significantly when the parameters of the  model are updated sequentially using training data of $\textsc{\textit{SC}}_{0.075}$, $\textsc{\textit{SC}}_{0.125}$ and $\textsc{\textit{SC}}_{0.2}$\footnote{Details of the dataset in Experimental Section}. Similar phenomenon is observed on other problems in the sequence such as $\textsc{\textit{SC}}_{0.075}$, etc.  It can be clearly concluded from Fig.~\ref{fig:intro_forget_sc} that GCNN~\cite{gasse2019exact} suffers from catastrophic forgetting in the lifelong learning  scenario. 
\looseness=-1

In the context of MIPs, a lifelong learning paradigm, due to its nature, promotes efficient learning since retraining from scratch is costly~\cite{parisi2019continual}, especially on industrial level MIPs where the number of variables and constraints are in orders of millions~\cite{nair2020solving}. Further, such a paradigm of \textit{continually gaining competencies} on different problems, offers opportunity to transfer the gained knowledge to unseen as well as previously seen problems. Owing to knowledge sharing across problems, these models cope better with low availability of training data. Learning from low-volume training data is important in the context of MIPs since generating training data (state-action pairs) itself is a computationally expensive process~\cite{nair2020solving}. Finally, using a single model that is updated regularly with various competencies, instead of maintaining multiple problem-specific models is desired as it is memory efficient with lower maintenance overhead. 

Motivated by the above listed benefits of lifelong learning, in this work we focus on the novel paradigm of \textit{Lifelong Learning to Branch in Mixed Integer Programs}. Our core contributions are as follows.

\noindent

\begin{itemize}
    \item \textbf{Problem Formulation:} We present the paradigm of Lifelong Learning to Branch in Mixed Integer Programs. To the best of our knowledge, we are the first to investigate this paradigm. 
    
    \item \textbf{Investigation of catastrophic forgetting:} We conduct an empirical investigation of the SoA method~\cite{gasse2019exact} and demonstrate that it suffers from catastrophic forgetting when it learns to branch on  different problems in succession.

    \item \textbf{Novel Algorithm:} We propose \namemodel, a \underline{Li}felong Learning method to solve \underline{M}ixed \underline{I}nteger \underline{P}rograms. \namemodel encodes the state-space of a problem through an \textit{edge-weighted, bipartite Graph Attention Network}. To mitigate catastrophic forgetting, \namemodel utilizes a novel combination of \textit{Knowledge Distillation} and \textit{Elastic Weight Consolidation} to shield the key parameters of previously learned problems from significant drift.

    \item \textbf{Experimental Evaluation:} We conduct empirical evaluation on a series of NP-hard problems with \textit{drifting} data distribution and evolving constraints. We establish that \namemodel is effective in learning to solve MIPs in a lifelong fashion and overcomes the problem of catastrophic forgetting. Further, it is also capable of transferring the gained knowledge effectively to NP-hard problems with very limited amount of training data.
\end{itemize}

%% file: 04_Method.tex
\section{Preliminaries}
\begin{defn}[Mixed Integer Program]
A mixed-integer linear program is an optimization problem of the form:
\begin{gather*}
\operatorname{minimize} \; \mathbf{c}^{\top} \mathbf{x} \\
\text{subject to}\; \mathbf{A} \mathbf{x} \leq \mathbf{b}  \\
\mathbf{l} \leq \mathbf{x} \leq \mathbf{u} \\
\mathbf{x} \in \mathbb{Z}^p \times \mathbb{R}^{n{-}p}
\end{gather*} 
where $n$ is the total number of variables, $p$ is the number of integral variables. $\mathbf{x} \in \mathbb{Z}^p \times \mathbb{R}^{n{-}p}$, $\mathbf{A} \in \mathbb{R}^{m \times n}$ is the constraint coefficient matrix, $\mathbf{b} \in \mathbb{R}^{m}$, right hand side constraint coefficient vector, $\mathbf{c} \in \mathbb{R}^{n}$ is the objective coefficient vector. Further, the variables $\mathbf{l}, \mathbf{u} \in  \mathbb{R}^n$ represent the lower and upper variable bound vectors. 

\end{defn}

MIPs are solved widely using Branch-and-Bound (B\&B) technique, which relaxes the integrality constraints and obtains a continuous linear program (LP). The LP is solved efficiently  using the simplex algorithm~\cite{achterberg2007constraint}. In case the relaxed solution is also integral and respects all the constraints, then it is also a solution to the problem (not necessarily optimal). Otherwise B\&B  decomposes the LP relaxation into two sub-problems, by splitting the feasible region based upon a variable that does not respect integrality constraints in the current LP solution $x^*$. Specifically, 
 
$$
x_{i} \leq\left\lfloor x_{i}^{\star}\right\rfloor \vee x_{i} \geq\left\lceil x_{i}^{\star}\right\rceil, \quad \exists i \leq p \mid x_{i}^{\star} \notin \mathbb{Z},
$$
 
The B\&B solving process repeatedly performs decomposition generating a search tree. The process stops if both the upper and lower bounds are equal or when the feasible regions do not decompose anymore, which is a certificate of infeasibility or  optimality. This B\&B procedure involves an extremely important step of selecting the fractional \textit{decision} variable to branch upon from the set of candidate variables $\mathcal{C}$. The chosen variable is used to partition the search space and has a significant impact on the size of the resulting search tree~\cite{achterberg2007constraint}. 

\looseness=-1
Among several heuristics available to choose the branching variable, \textit{strong branching} is widely known to produce the smallest B\&B trees. It calculates the expected bound improvement for each candidate variable before performing branching. Although it produces the smallest B\&B trees, strong branching  requires computing the solution of two LPs for each candidate variable. The cost of finding the best variable is prohibitively high and hence strong branching is not used in practice. 

 In the B\&B setup, the MILP solver is considered to be the \textit{environment}, and the brancher the \textit{agent}~\cite{khalil2016learning,gasse2019exact}. At the $t^{th}$ decision step, the solver is in  state ${s_t}$, which comprises the B\&B tree with all past branching decisions~\cite{gasse2019exact}, the best integer solution found so far, the LP solution of each node, the currently focused leaf node, as well as any other solver statistics (for example, the number of times every primal heuristic has been called). In the context of strong branching, at a given state $s_t$, let $a_t$ be the  variable chosen by strong branching among the set of all candidate variables $C$. Based upon the above discussion, we now define the problem of \textit{Learning to Branch}.

\looseness=-1
\begin{prob}[Learning to Branch]
\looseness=-1 For a B\&B tree, at the $t^{th}$ decision step of the solver, let the solver be in state ${s_t}$ and a decision to choose a variable to branch is to be made from a set of candidate variables. Given a collection of  state-action $({s_t}, {a_t} )$ pairs obtained from running strong branching, the goal is to learn a scoring function $f$ parameterized by $\theta$ that imitates branching decisions made by the strong branching expert. 

\end{prob}

\noindent
Since in our setup we aim to learn on multiple problems, therefore in this context, we refer to each problem as a task. 

\begin{figure*}[t!]
\centering
\includegraphics[scale=0.45]{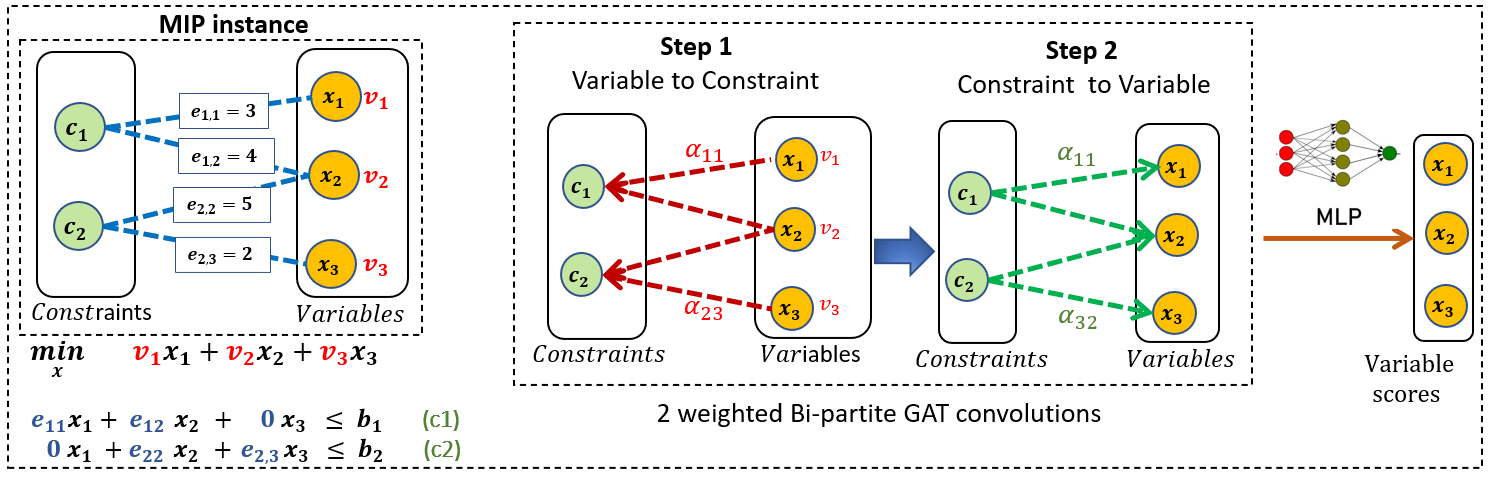}
\caption{Bipartite graph representation of an MIP  with $n=3$ variables and $m=2$ constraints. The bipartite graph is encoded via 2 half-aggregations of Bipartite GAT .}
\label{fig:architecture}
 \vspace{-0.20in}
\end{figure*}
\begin{prob}[Lifelong Learning to Branch]
Given a sequence of tasks $\mathcal{T} =[ \mathcal{T}_1,\cdots, \mathcal{T}_T ]$ of length $T$, we aim to  update the parameters of the model sequentially over time such that at the  $i^{th}$ task, when parameter $\theta_{i-1}$ is updated to  $\theta_i$ by training using the  instances of the task $\mathcal{T}_i$,  the model avoids catastrophic forgetting on tasks $\mathcal{T}_j$ for $j{<}i $. Specifically, the increase in running time on problems $\mathcal{T}_j \;  \forall  j{<}i $ using the updated model $\theta_i$ should be reasonably low. Additionally, the performance on newly learned tasks should also not be hindered significantly. 
 \label{def:lifelongbranch}
\end{prob}

\section{\namemodel: Our Proposed Metholodogy}

In this section we describe our proposed method \namemodel. Fig.~\ref{fig:architecture} presents the overview of \namemodel. We first convert a given MIP instance to a bipartite graph. Next, we describe a method to encode the variables and constraints of MIP using an \textit{edge-weighted, bipartite graph attention network (GAT)}. Finally, we describe the procedure of learning the parameters of the model in a lifelong  fashion by avoiding catastrophic forgetting. With the outline being set, we next discuss each of these components in detail.

\subsection{MIP Representation: State Encoding\\}

Similar to \citet{gasse2019exact}, to encode the state  $\mathbf{s_t}$  of the B\&B tree at timestep $t$,  we use a bipartite graph representation $G = (\mathbf{V}, \mathbf{E}, \mathbf{C})$ . One side of  the graph containing $n$ nodes represent the $n$ variables and the other side consisting of $m$ nodes represent the $m$ constraints. There exists an edge between   $j^{th}$ variable node  and $i^{th}$ constraint node if  $j^{th}$ variable appears in the $i^{th}$ constraint. The weight of an edge $e_{i,j}$ corresponds to the value of the coefficient of the variable $v_j$ in the constraint $c_i$.  We use $\mathbf{V} \in R^{n\times d_1}$ to represent variable features, $\mathbf{E} \in R^{n  \times m \times  1}$  for edge features, and $\mathbf{C} \in R^{m\times d_2}$ to represent the constraint features. For each of the node in the graph we use the raw solver specific input features of \citet{gasse2019exact}, which can be found in the appendix A.2. 
  Fig.~\ref{fig:architecture} shows an example of encoding an MIP to a bipartite graph.

\subsection{Policy Parameterization: \textit{Edge Weighted Bipartite GAT}}

Observing the weighted and bipartite nature of the graph, it is natural to parameterize the branching variable policy $f_\theta(a|s_t)$ using an \textit{edge-weighted} bipartite GAT. Specifically, for each node in the graph, the attention layer learns to weigh each of the node's neighbors differently based upon its importance~\cite{vaswani2017attention}. Since, our graph is bipartite, we perform two levels of message passing through the GAT. Specifically, first we pass message from the variable side  to the constraint side to obtain rich representation of the constraint nodes as follows:
\vspace{-0.05in}
\begin{equation}
\label{eq:conv1}
\mathbf{c}_{i}=\alpha_{i, i} {\theta^C} \mathbf{c}_{i}+\sum_{j \in \mathcal{N}(i)} \alpha_{i, j} \theta^C \mathbf{v}_{j}    
\end{equation}

Here, $\mathbf{c}_{i}$ and $\mathbf{v}_{j}$ refer to the embeddings of  $i^{th}$ constraint and  $j^{th}$ variable respectively. $\mathcal{N}(i)$ refers to the neighbors of $i^{th}$ node. $\theta^C$ refers to MLP associated to constrained side aggregation. $\alpha$ represents the attention coefficient(defined later). Next, we perform message passing from constraint side to variable side. This allows us to generate richer  representations for each of the variables nodes. 
\begin{equation}
    \label{eq:conv2}
    \mathbf{v}_{j}=\alpha_{j, j} {\theta}^V \mathbf{v}_{j}+\sum_{i \in \mathcal{N}(j)} \alpha_{j, i} \theta^V \mathbf{c}_{i}
\end{equation}

$\theta^V$ refers to weights associated to the variable side aggregation. 
The attention coefficient  $\alpha$ is computed as below:

$$
\alpha_{{i}, {j}}=\frac{\exp \left(\rho\left({\mathbf{(a^C)}}^{\mathrm{T}}\left[\theta^{C} {\mathbf{c}_{i}} \| \theta^{C} {\mathbf{v}_{j}}  \| \theta^C_{e} \mathbf{{e}}_{i,j} \right]\right)\right)}{\sum_{{k} \in \mathcal{N}_{(i)}} \exp \left(\rho\left({(\mathbf{a^C})}^{\mathrm{T}}\left[\theta^{C} {\mathbf{c}_{i}}  \| \theta^{C} {\mathbf{v}_{k}} \|  \theta^C_{e} \mathbf{{e}}_{i,k} \right]\right)\right)} $$

\looseness=-1
The above attention mechanism is parameterized by the weight vector $\mathbf{a}^C$.    $\theta^C_{e}$ on the constraint side aggregation  refers to an MLP associated with the edge features. $\rho$ refers to the activation function\footnote{We use LeakyReLU with negative slope $=0.2$}.  $\alpha_{{j}, {i}}$ is defined analogously to $\alpha_{{i}, {j}}$ where ${C}$ is swapped with ${V}$ and the $i^{th}$ and $j^{th}$ nodes are interchanged.
 
 After the two half-aggregations of eq.~\ref{eq:conv1} and ~\ref{eq:conv2}, we obtain the final representation of the candidate variable nodes.  The final representation of each candidate variable node is passed through a softmax layer to obtain a probability distribution over the variables represented by $f_{\theta}\left({a} \mid {s_t}\right)$. 
 Further, to stabilize the training procedure of the bipartite-GAT, we  use attention mechanism with multiple heads,  details of which are present in  Appendix Sec.A.3.  
The detailed architecture is present in  Fig.~\ref{fig:architecture}.

\vspace{-0.05in}
\subsection{Imitation Loss}
Since strong branching is a powerful heuristic in reducing the size of the tree, we train the parameters $\theta$ by imitation learning of the strong branching rule. 

We first collect a set of strong branching state-action pairs  $\mathcal{D}=\left\{\left({s}_{k}, {a}_{k}^{\star}\right)\right\}_{k=1}^{N}$, where $N$ is the number of branching samples collected.  Then through imitation learning, we optimize the parameters $\theta$  using the following imitation loss function:
\vspace{-0.05in}
\begin{equation}
\label{eq:Loss_c}
\mathcal{L}(\theta)=-\frac{1}{N} \sum_{\left({s}, a^{*}\right) \in \mathcal{D}} \log f_{\theta}\left({a}^{*} \mid {s}\right)
\end{equation}

The optimization of the above objective encourages the neural model to predict the  variable for branching which strong-branching would have chosen.

\subsection{Life-Long Learning to Branch}
Until now we discussed how to learn the parameters $\theta$ of the model for a given task. In this section we discuss how to learn to branch on MIPs in a lifelong fashion. As discussed earlier in Def.~\ref{def:lifelongbranch}, we have a set of $T$ problems appearing in sequence $\mathcal{T} =[ \mathcal{T}_1,\cdots, \mathcal{T}_T ]$, and our goal is to learn the parameters sequentially over time where the training data $\mathcal{D}_i$ for each task $\mathcal{T}_i$ also appears sequentially. A na\"ive solution is to update the parameters of the model sequentially as new tasks arrive. However, as we already observed in Fig.~\ref{fig:intro_forget_sc}, if the neural model is updated in this fashion, it suffers from catastrophic forgetting on the earlier learned tasks. Hence, our goal is to update the parameters of the model on new tasks while preserving the knowledge gained on previous tasks to avoid catastrophic forgetting.

\begin{figure}[b]
\centering
\includegraphics[scale=0.21]{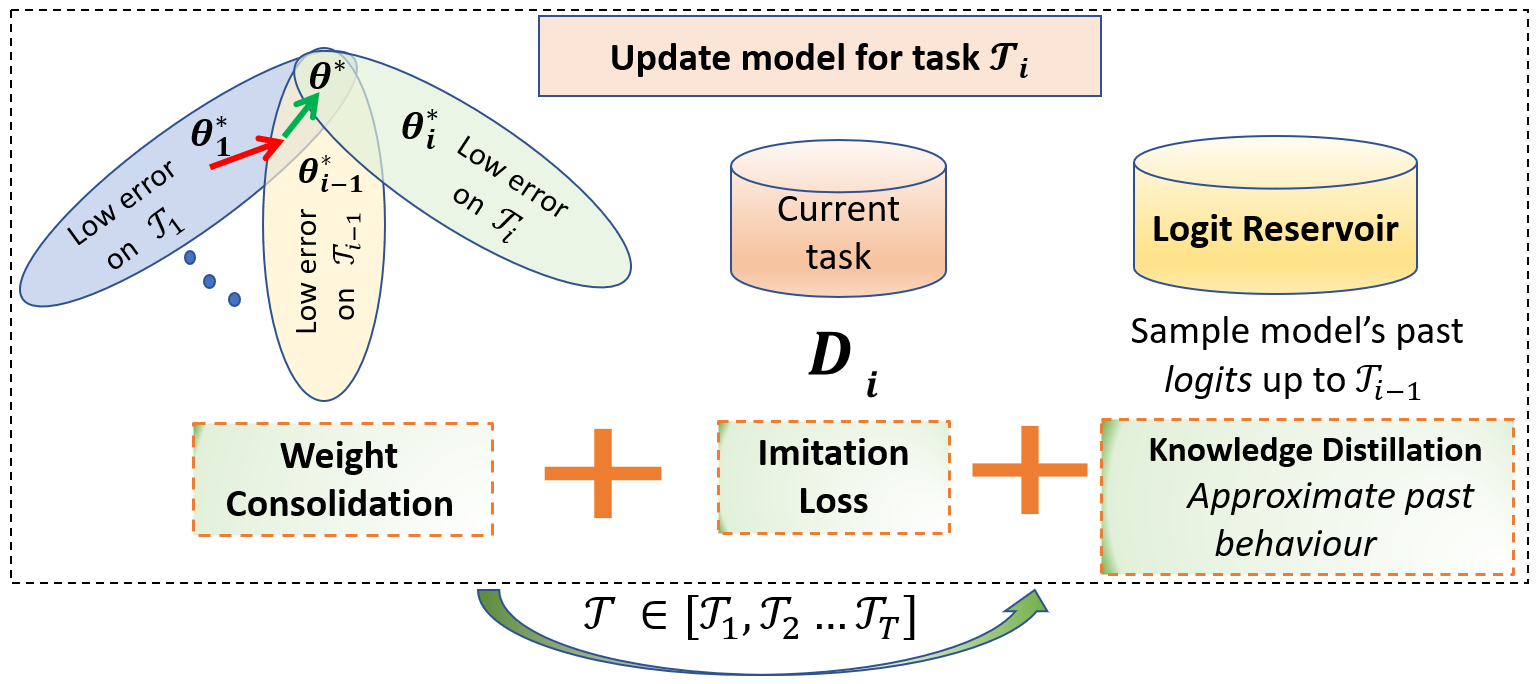}
\caption{Architecture diagram representing the update mechanism of \namemodel at the $i{^{th}}$ step of sequence $\mathcal{T}$.}
\label{fig:architecture_life}
\end{figure}

One way to consolidate past knowledge is to replay the training data of the past tasks. However, as the number of tasks increase, it becomes computationally expensive. Further, another option of storing only a small set  of labeled samples and replaying them is prone to over-fitting~\cite{wang2020streaming}. Hence, inspired by recent works on continual learning~\cite{buzzega2020dark, wang2020streaming}, to tackle the problem of catastrophic forgetting in \textit{lifelong learning to branch} we take the following two perspectives. First idea is to approximate the knowledge gained by the model in the past via \textit{distillation} of model's past behavior when learning new tasks. Second, we optimize the parameters of the model in a constrained way in order to prevent significant drift on the parameters important for previously learned tasks. Fig.~\ref{fig:architecture_life}. visually describes the process. Now, we discuss both the perspectives below in detail.

\subsubsection{Mimicking model's past behavior through  Knowledge Distillation:}

In order to maintain past learned patterns during lifelong learning, our goal is to search for model parameters that fit well on the current task and also approximate the optimal behavior of the model on the older tasks. 
Towards this we aim to encourage the model to mimic its original(past) output \textit{logits} for a small number of samples of the past tasks. To accomplish this we apply  \textit{Knowledge distillation}(KD)~\cite{buzzega2020dark} approach to enforce the neural network to generate similar logits that the model produced for these samples in the past during optimizing of the task to which the related sample belonged to. Mathematically, 

\vspace{-0.20in}
\begin{equation}
    \label{eq:dataKLD}
    \mathcal{L}_{KL} =  \mathbb{E}_{(s,z) \sim M }\left[D_{K L}\left(z \| f_{\theta}(s)\right)\right]
\end{equation}

Here $z=f_{\theta_j*}(s)$ refers to the logits $z$ of sample $s$  and  $\theta^{*}_j$ refers to the set of optimal parameters of task $\mathcal{T}_j$. These $(s,z)$ pairs are stored in a fixed-size buffer $M$. Specifically, for $s \in M$, $f_{\theta_{j}^{*}}(s)$ is preserved where $s$ is a training sample from task $\mathcal{T}_j$. When the lifelong learning model is at step $i$ of the sequence $\mathcal{T}$, $M$ consists of samples of past experiences(logits) for tasks seen till step  ${i-1}$. 
Further, since we do not have any prior information of how many tasks we will observe, we use \textit{reservoir sampling} to preserve samples for Knowledge-Distillation. Reservoir sampling ensures that samples from all tasks are stored with equal probability in the buffer without knowing the number of tasks/samples in the stream in advance\cite{buzzega2020dark}.

\subsubsection{Preservation of model's important parameters:}
As we store only a small set of logits in our memory buffer instead of the entire training data of past tasks, it is prone to over-fitting. Although, over-fitting can be tackled to an extent by L2 regularizers, the restriction imposed by L2 regularizers by constraining the entire network through a fixed coefficient is too severe and might prevent learning of the new tasks itself. Inspired from recent works~\cite{kirkpatrick2017overcoming,wang2020streaming}, to counter this problem, we aim to learn to adjust the magnitude of the parameter updates on certain model weights based on how important they are to the  previously learned tasks. To accomplish this we apply \textit{Elastic Weight Consolidation   (EWC)}~\cite{wang2020streaming,kirkpatrick2017overcoming}. Specifically, after the training on a task $\mathcal{T}_j$ is complete, we compute the importance of each parameter $w$ on the task $\mathcal{T}_j$ as follows:

\begin{center}

$
\Omega_{j}^{w}=\mathbb{E}_{(s,{a^*}) \sim D_{j}}\left[\left(\frac{\delta \mathcal{L}(s,{a^*})}{\delta \theta_{j}^{w}}\right)^{2}\right]
 $       
 
\end{center}

 $\mathcal{L}(s,{a^*})$ refers to the loss on sample $s$ with ground-truth $a^*$.  The term $\frac{\delta \mathcal{L}(s,{a^*})}{\delta \theta_{j}^{w}}$ calculates the gradient of the loss  with respect to the parameter $w$ learned on task $\mathcal{T}_j$. $\Omega_{j}^{w}$ captures the importance of weight $w$ to task $\mathcal{T}_j$. We note that $D_j$ is no more required during future tasks once the computation of $\Omega$ for task $\mathcal{T}_j$ is complete.

Now, when a new task $\mathcal{T}_i$ arrives, we apply the above regularization (penalize) to prevent large amount of drift on parameters important for earlier learned tasks. Here, the weights of the regularization are obtained from $\Omega$. We accomplish the regularization by the below loss function.
\vspace{-0.05in}

\begin{equation}
\label{eq:model}
\mathcal{L}_{importance} = \sum_{j=1}^{i-1} \sum_{w}  \Omega_{j}^{w}\left(\theta_{i}^{w}-\theta_{j^*}^{w}\right)^{2}
\end{equation}


The above term is a quadratic penalty term on the difference between the parameters for the new and the old tasks. $\Omega$ consists of diagonal weighing proportional to the diagonal of the Fisher information metric over the old parameters on the old tasks~\cite{liu2021overcoming}. $\theta_{j^*}$ refers to optimal parameters of task $\mathcal{T}_j$ . 
When updating parameters of the model to learn to branch on a new task $\mathcal{T}_i$, the above penalization will encourage the \textit{important model} parameters to be close to the parameters obtained for earlier learned tasks $\mathcal{T}_1, \mathcal{T}_2 \cdots \mathcal{T}_{i-1}$. Fig.~\ref{fig:architecture} b) summarizes this concept visually through overlapping optimal parameter spaces.

\subsubsection{Lifelong learning optimization objective:}
Finally, combining the loss functions of eqs.~\ref{eq:Loss_c},~\ref{eq:dataKLD} and ~\ref{eq:model}, we obtain the optimization objective at the $i^{th}$  step as follows

\vspace{-0.05in}

{
\begin{equation}
\begin{aligned}
\nonumber
 \mathcal{L}_{lifelong} =&  \sum_{\left({s}, {a}^{*}\right) \in \mathcal{D}_{i}} \log f_{\theta}\left({a}^{*} \mid {s}\right)
 \\& + \alpha \mathbb{E}_{(s, z) \sim \mathcal{M}} [ D_{K L}\left(z \| f_{\theta}(s) \right]  
 \\& + \beta \sum_{j=1}^{i-1} \sum_{w}  \Omega_{j}^{w}\left(\theta_{i}^{w}-\theta^{w}_{j^*}\right)^{2}
\label{eq:overall_loss}
\end{aligned}
\end{equation}}
\looseness=-1
The above equation while learning new tasks, consolidates past information in order to maintain stability of parameters important for previously learned tasks. $\alpha$ controls the weight corresponding to mimicking past logits and $\beta$ controls scale of the weight consolidation regularizer. The detailed steps of training a sequence of tasks in lifelong fashion through \namemodel are described in Algorithm 1 in Appendix A.4.

%% file: 05_Experiments.tex
\section{Experiments}
In this section we measure the effectiveness of our proposed approach ~\namemodel and establish:
\begin{itemize}
    \item \textbf{Minimal forgetting:} \namemodel is capable of lifelong learning on NP-hard problems with drifting data distributions and avoids catastrophic forgetting on previously learned problems.
    
    \item \textbf{No hindrance in learning future tasks:} Despite adding constraints to prevent significant updates to the model, \namemodel does not hinder learning on new tasks.
    
    \item \textbf{Transfer to Low data regime:} We compare the performance of \namemodel to transfer on a low-training data regime task which is similar to a task \namemodel learned in the past. \namemodel effectively transfers its previously gained and \textit{unforgotten knowledge} to the unseen task.
    
    \item \textbf{Efficient learning through Bipartite GAT:} Attributed to rich representations learned through the attention mechanism, \namemodel reduces solving time on instances when compared to \GCNN~\cite{gasse2019exact}.
    
\end{itemize}

The codebase can be found on \url{https://github.com/idea-iitd/LiMIP} .

\subsection{Datasets}
We use the following datasets to evaluate the performance of our method against different baselines.

\looseness=-1
\textbf{Set Cover}: \looseness=-1We consider the Set Cover problem of \citet{balas1980set}. Let $p$ be the probability of an item belonging to a set in the Set Cover(\textsc{$SC$}) problem. \textsc{$SC_{p}$} refers to Set Cover problem with set-item probability $p$. To simulate lifelong learning setup, we generate multiple Set Cover problems datasets each with a different probability, i.e., $\mathcal{T} = \left[ \textsc{$SC_{0.05}$}, \textsc{$SC_{0.075}$}, \textsc{$SC_{0.1}$}, \textsc{$SC_{0.125}$}, \textsc{$SC_{0.15}$}, \textsc{$SC_{0.2}$} \right]$. In all instances we set number of rows to $700$ and number of columns to $800$.

\textbf{Independent Set }: We consider the Maximum Independent Set (MIS) problem on the Barabási-Alberta graph~\cite{albert2002statistical} generated with different sizes and affinities.  \textsc{$IS_{{A},{S}}$} denotes as instance where $A$ is the affinity and $S$ is the size of the graph. To simulate lifelong learning setup, we generate independent set problem datasets  with  different sizes and affinities as  $\mathcal{T} = \left[ \textsc{$IS_{{4},{750}}$} , \textsc{$IS_{{4},{500}}$}, \textsc{$IS_{{4},{450}}$}, \textsc{$IS_{{5},{450}}$}, \textsc{$IS_{{5},{400}}$}, \textsc{$IS_{{5},{350}}$}  \right]$.

\textbf{Facility Location with constraints:}
We consider the Facility Location problem~\cite{gasse2019exact} and to simulate lifelong learning scenario, we use facility capacities and customer demands sampled from drifting distributions over time. The details can be found in Appendix A.5. 


These datasets are challenging for state-of-the-art solvers, and  also representative of the types of integer programming problems encountered in practice.

\looseness=-1
\subsection{Experimental Setup and Parameters}
\looseness=-1We use SCIP~\cite{gamrath2020scip} as the backend solver, with a time limit of 45 minutes. We use a system running on Intel Xeon 6248 processor with 96 cores and 1 NVIDIA A100 GPU with 40GB memory for our experiments. Similar to existing works~\cite{gasse2019exact}, we enable cutting planes at the root node and deactivate solver restarts. We keep all other  SCIP parameters to default.  We use attention mechanism with $2$ heads. We set the default buffer size to $500$. For details of all parameters and system settings, we refer to App A.6. 

\noindent
\textbf{Training data generation:}
For each task, we generate 150,000  branching samples extracted using 10,000 generated instances for training and 30000 validation/test samples generated using 2000 instances.

\subsubsection{Metrics:} We perform evaluation on 20 different test instances using 5 different SCIP seeds.  We report the standard benchmark metric for MILP benchmarking, i.e., the geometric mean of the  \textbf{running time} of the solver. Additionally, we report the\textbf{ hardware independent node count} (in Appendix A.8). 
 We  compute the average per-instance standard deviation so a value $X \pm s  \% $ means it took $X$ secs to solve an instance and while solving one of those instances the time varied by $s$ on an average.

\subsubsection{Baselines:}
We compare our work \namemodel with the state-of-the-art method for learning to branch \textbf{GCNN}~\cite{gasse2019exact}.   We skip comparison with \citet{zarpellon2020parameterizing} since it approximates the weaker reliability pseudocost branching heuristic, which has been shown to have an inferior performance in terms of running time~(See Appendix A.10). 
 Further, we skip comparison with ~\citet{gupta2020hybrid} since its focus is on developing CPU based version of learning to branch, which is out of scope of our work.  For the sake of completion, we compare with the default SCIP Solver and strong branching in App. A.9. 

In the context of lifelong learning, we compare with \textbf{(1)} Fine-tuning (FT) i.e directly updating the model on new tasks as they arrive, \textbf{(2)} Experience Replay (ER) and \textbf{(3)} Elastic Weight Consolidation (EWC)~\cite{kirkpatrick2017overcoming}. Details of baseline are present in App. A.7. 

\begin{figure}[t]
\hspace{-0.15in}
\includegraphics[scale=0.44]{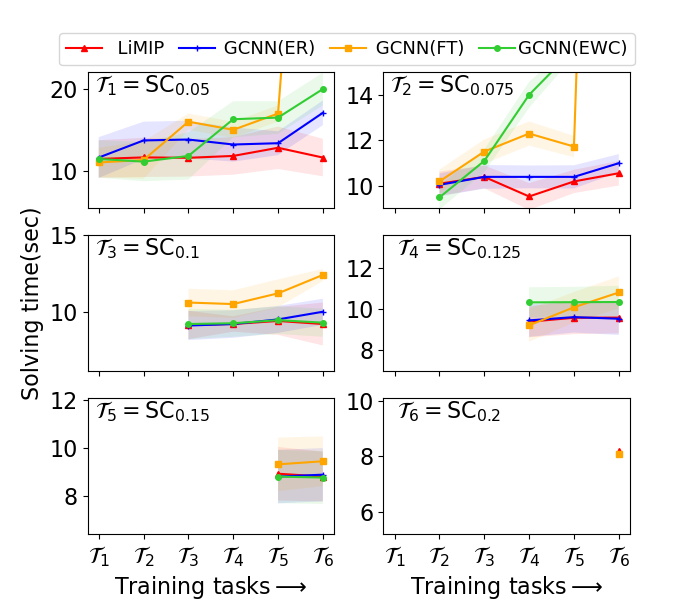}
\caption{{Testing on Set Cover in lifelong scenario}: Evolution of solving time  for each task when different methods are updated on each task sequentially. Different \textit{evaluation tasks} are shown in \textit{different subplots}. The x-axis denotes the sequence of training tasks and the y-axis denotes the geometric mean of solving time for \textit{test instances} of each task in the sequence.  The shaded area refers to standard dev.}
\label{fig:eval_sc_continue}
\end{figure}
\subsection{Evaluation in Lifelong Learning Scenario}

\begin{figure}[b!]
\centering
\hspace{-0.3in}
\includegraphics[scale=0.41]{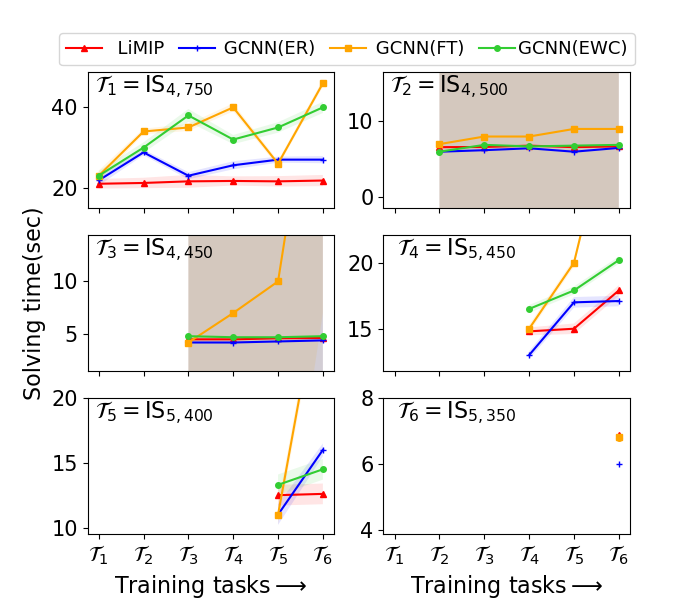}

\caption{ {Testing on Independent Set in lifelong scenario}: Evolution of solving time  for each task when different methods are updated on each task sequentially.}
\label{fig:eval_fac_continue}
\end{figure}

\textbf{Evaluating forgetting:} In  Fig.~\ref{fig:eval_sc_continue} and ~\ref{fig:eval_fac_continue}, in each subplot we study the performance of different methods on  test instances of each  dataset in the lifelong sequence $\mathcal{T}$. Specifically, each subplot in these figures refer to a test task and the x-axis shows the sequence of training tasks in the lifelong setup.  We observe that as training progresses on different problems in the lifelong setup, the performance of \GCNN(FT)~\cite{gasse2019exact} on old problem deteriorates significantly. The older the task, the worse is the deterioration. For example, since $SC_{0.05}$ (top, left in Fig.~\ref{fig:eval_sc_continue}) is trained first, it witnesses maximum increase in time across algorithms except the proposed \namemodel. This clearly shows that \GCNN suffers from catastrophic forgetting when its parameters are updated on new tasks.  In sharp contrast, \namemodel is able to maintain the learned patterns of past tasks when learning new tasks. This is attributed to the \textit{knowledge distillation} loss, which helps in promoting the model to  mimic past behavior and the \textit{weight consolidation penalty} term which prevents significant drift on important parameters. 
 
\textbf{Low hindrance on future tasks:} 
\looseness=-1 While avoiding catastrophic forgetting is one of the aims of lifelong learning, it should not be at the cost of learning new tasks.  From  Fig.~\ref{fig:eval_sc_continue} and Fig.~\ref{fig:eval_fac_continue}, we observe that, while \namemodel does not forget the knowledge it gained in the past, still it does  it without impacting  future tasks. We can clearly see that on future tasks too \namemodel obtains superior performance compared to existing baselines. In App A.11 
we also compare with GAT (FT), GAT (EWC) and GAT (ER). The results on Facility Location Problem are present in App A.9. 

\subsection{Transferability on Low Data Regimes}
\looseness=-1
As a model learns on more tasks and gains several competencies, it can be utilized as a \textit{weight initializer} to learn on an unseen task with low availability of training data, especially in the case of strong branching where obtaining training data is costly. To test the performance of \namemodel on transfer learning, we create  a new dataset $SC_{0.047}$  with only $300$ branching samples for training. 
This number is extremely low compared to  other datasets where number of samples $=150000$. We perform lifelong learning on  $\mathcal{T}= \left[ \textsc{$SC_{0.05}$}, \textsc{$SC_{0.075}$}, \textsc{$SC_{0.1}$}, \textsc{$SC_{0.125}$}, \textsc{$SC_{0.15}$}, \textsc{$SC_{0.2}$} \right]$ and then fine-tune on \textsc{$SC_{0.047}$}.  We compare it with training a model from scratch on the $300$ branching samples of \textsc{$SC_{0.047}$}. In Table~\ref{tab:transfer}, we study the performance gain obtained using a fine-tuned lifelong learned model vs. a model learned from scratch $(\GCNN \; \textsc{$SC_{0.047}$}$). We observe that on an average the running time of the \namemodel method fine-tuned on $\textsc{$SC_{0.047}$}$ model is significantly better than the model trained from scratch.

\begin{table}[t]
    \centering
        \begin{tabular}{p{4.2cm}p{1.3cm}p{1.3cm}p{1.3cm}}
\toprule
        \textbf{Model} & \textbf{Time(sec)} & \textbf{\# Nodes}  \\
        \midrule
        \GCNN   \textsc{$SC_{0.047}$}  &   {$16.01$} &  {$607$} \\ 
       { Bipartite GAT  $SC_{0.047}$} &  {$17.06$} &  {$702$}\\
       { Bipartite  GAT(FT) + $SC_{0.047}$} &  {$16.07$} &  {$702$}\\
        \namemodel   &   {$15.05$} &  {$602$} \\
\namemodel + ($SC_{0.047}$)  &   $\mathbf{14.05}$ &  {$\mathbf{441}$} \\
        \bottomrule
        \end{tabular}
    \caption{{Transferability performance}: Test performance comparison of fine-tuned model against model trained from scratch on the $SC_{0.047}$ dataset. \namemodel and Bipartite GAT (FT) were trained sequentially on  $\left[ \textsc{$SC_{0.05}$}, \textsc{$SC_{0.075}$}, \textsc{$SC_{0.1}$}, \textsc{$SC_{0.125}$}, \textsc{$SC_{0.15}$}, \textsc{$SC_{0.2}$} \right]$ and then fine-tuned on  $\textsc{$SC_{0.047}$}$ }
    \label{tab:transfer}
\end{table}

\subsection{Ablation Studies}
\subsubsection{GCN vs Bipartite GAT:}
In this section we study the impact of using our Bipartite GAT compared to mean pool based \GCNN. In Table~\ref{tab:res:GATGCN} we observe that Bipartite GAT improves over GCNN by a small margin in terms of both running time and number of tree nodes.

\setlength{\tabcolsep}{2pt} 
\begin{table}[h!]

\centering
\small
\begin{tabular}{llcc}

\toprule
\textbf{Dataset} & \textbf{Method} & \textbf{Time(sec)} & \textbf{Nodes}\\ 
\midrule
\multirow{2}{*}{\textsc{$SC_{0.2}$} } & Bipartite GAT & \boldmath{$7.9 \pm 2.56$} & \boldmath{$87.2 \pm15.31$}\\
 & \GCNN & $8.21 \pm2.45$ & {$91.2 \pm13.61$} \\
\cmidrule{1-4}

\multirow{2}{*}{\textsc{$IS_{4},{750}$} } & Bipartite GAT & \boldmath{$22.25 \pm1.24$} & \boldmath{$555.6 \pm6.02$}\\
 & \GCNN & {$25.73 \pm 1.32$} & {$672.2 \pm6.2$} \\
\cmidrule{1-4}

\multirow{2}{*}{\textsc{$FC_{(40,50),(5,10)}$} } & Bipartite GAT & \boldmath{$33.94 \pm1.11$} & \boldmath{$246.10 \pm4.02$}\\
 & \GCNN & {$35.14 \pm 1.30$} & {$248.20 \pm4.70$} \\
\cmidrule{1-4}

\bottomrule

\end{tabular}
\caption{Performance comparison between Bipartite GAT  and \GCNN encoding.}
\label{tab:res:GATGCN} 
\end{table}

\subsubsection{Impact of  regularizer and buffer size on lifelong learning:} In App A.11 
 we study the impact of weight regularization and buffer size. 


%% file: 06_Conclusion.tex
\section{Conclusion}
\vspace{-0.05in}
Learning-to-Branch techniques have shown significant success in reducing the solving time of Mixed Integer Programs. Although, significant progress has been made, the paradigm of \textit{learning to branch} in a \textit{lifelong fashion} was unexplored. In this work we first examined the behavior of existing techniques in the lifelong learning scenario and discovered that they suffer from catastrophic forgetting. To tackle this problem, in this work we study the hitherto unexplored paradigm of \textit{Lifelong Learning to Branch} on Mixed Integer Programs. We propose a method \namemodel powered by a Bipartite GAT to encode MIP instances. Further, to mitigate catastrophic forgetting, we apply \textit{knowledge distillation} and \textit{elastic weight consolidation} to shield key parameters from drifting and thereby retaining efficacy. Through extensive experiments on multiple NP-hard problems, we established that \namemodel is able to mitigate forgetting significantly better compared to existing baselines when confronted with lifelong learning. Additionally, the proposed method does not hinder the performance on future learning tasks too. 

%% file: 07_ack.tex
\section*{Acknowledgements}
We acknowledge  Qualcomm for supporting this research. We acknowledge Google for supporting Sahil Manchanda for this travel.

%% file: Appendix.tex
\clearpage
\section{Appendix}



\subsection{Related Work}
\label{app:relatedwork}
\input{03_RelatedWork}
\subsection{Solver Features}
\label{app:solver_features}
In Table ~\ref{tab:features} we present the input features used. These features are based upon ~\citet{gasse2019exact}.

\renewcommand{\tabcolsep}{3pt}
\begin{table*}[h!]
\centering

\scalebox{1}{
{\scriptsize
\begin{tabular}{c|c|c}
    \textbf{Tensor} &\textbf{Feature} &  \textbf{Description}  \\
    
    \midrule
    \multirow{5}{2em}{$\mathbf{C}$} &  obj\_cos\_sim &Cosine similarity with objective. \\
    & bias & Bias value, normalized with constraint coefficients. \\
    & dualsol\_val & Dual solution value, normalized.\\
    & is\_tight &Tightness indicator in LP solution\\
    & age & LP age, normalized with total number of LPs \\
     \hline
    \hline
      \multirow{1}{2em}{$\mathbf{E}$} &  Coef & Constraint coefficient, normalized per constraint\\
    \hline
          \multirow{13}{2em}{$\mathbf{V}$} &  type & Type (binary, integer, impl. integer, continuous) as a one-hot encoding\\
  & coef & Objective coefficient, normalized \\
  & has\_lb & Lower bound indicator.\\
  & has\_ub & Upper bound indicator.\\
  & sol\_is\_at\_lb & Solution value equals lower bound. \\
  & sol\_is\_at\_ub & Solution value equals upper bound \\
  & sol\_frac & Solution value fractionality \\
  & basis\_status & Simplex basis status (lower, basic, upper, zero) as a one-hot encoding \\
  & reduced\_cost & Reduced cost, normalized \\
  & age & LP age, normalized \\
  & sol\_val &  Solution value. \\
  & inc\_val & Value in incumbent \\
  & avg\_inc\_val& Average value in incumbents\\



 
    \hline   

    \bottomrule
\end{tabular}}
}
\caption{\label{tab:features} Description of the constraint, edge and variable features in our bipartite state representation $s_t$}
\end{table*}

\subsection{Multi-head attention}
\label{app:MHA}
The multi-head attention equation is based upon ~\citet{vaswani2017attention}. Specifically, for $K$ heads
\begin{equation}
\nonumber
\vec{h}_{i}^{\prime}=\|_{k=1}^{K} \rho\left(\sum_{j \in \mathcal{N}{(i)}} \alpha_{i j}^{k} \mathbf{W}^{k} \vec{h}_{j}\right)
\end{equation}

where $h_i$ can be replaced by $c_i$ or $v_i$ and $\mathbf{W}$ be replaced with $\theta^C$ or $\theta^V$ respectively.

\subsection{Pseudocode}
\label{app:algo}
\begin{algorithm}
\caption{Lifelong Learning of parameters $\theta$ at $i^{th}$ step. }
\label{alg:lifelong_train}
\begin{algorithmic}[1]
 	 \REQUIRE  $\theta$, $D_i$, $M$, $\alpha,\beta$ 
	    \WHILE{not converged}
	          \STATE  $\mathcal{L}_i{=}-\frac{1}{N} \sum_{\left({s}, a^{*}\right) \in \mathcal{D}_i} \log f_{\theta}\left({a}^{*} \mid {s}\right)$
	          \COMMENT{Compute loss for current task $\mathcal{T}_i$}
	        
	        \STATE $\mathcal{L}_{KL}{=} \mathbb{E}_{(s,z) \sim M }\left[D_{K L}\left(z \| f_{\theta}(s)\right)\right]$ \COMMENT{Compute KL divergence  for samples $\in M$ }

	        \STATE  $\mathcal{L}_{lifelong}{=}\sum_{\left({s}, {a}^{*}\right) \in \mathcal{D}_{i}} \log f_{\theta}\left({a}^{*} \mid {s}\right) \newline \null \quad \quad \quad \quad \quad \quad {+}  \alpha \mathbb{E}_{(s, z) \sim \mathcal{M}} \left[D_{K L}\left(z\| f_{\theta}(s)\right) \right] \newline \null  \quad \quad \quad \quad \quad \quad  {+} \beta \sum_{j=1}^{i-1} \sum_{w} \Omega_{j}^{w}\left(\theta_{i}^{w}-\theta^{w}_{j^*}\right)^{2}$
 
          $\;$  Update parameters $\theta$ using loss $\mathcal{L}_{lifelong}$  $\;\;\;$
            
    \ENDWHILE

	        \STATE $\Omega_{i}^{w}{=}\mathbb{E}_{(s,{a^*}) \sim D_{i}}\left[\left(\frac{\delta \mathcal{L}(s,{a^*})}{\delta \theta_{i}^{w}}\right)^{2}\right]$
	        \COMMENT{Compute importance of weights at end of the task}
	        
	        \STATE{Update memory $M$ using reservoir sampling}
	        
	        
\end{algorithmic}
\end{algorithm}

\subsection{Data}
\label{app:data}
\textbf{Facility location:}
This is a realistic scenario where customer demands and facility capacities keep changing over time. 
For facility location, as discussed the main paper, a certain task of facility location problem is defined as $FC_{{(Clow, Chigh)},{(Dlow, Dhigh)}}$.  $Clow$, $Chigh$ refers to the lower and upper limit( respectively) of the facility capacity. $Dlow$, $Dhigh$ refers to the lower and upper limit(respectively)  of the customer demand . To generate an instance of the facility location problem, we  require, for each facility, sampling a capacity uniformly at random from $[Clow, Chigh]$ and for each customer a demand sampled uniformly at random from $[Dlow, Dhigh]$. Using these values, an instance is constructed.  In addition to evolving customer demand and facility supply distribution, we also simulate the setting of adding new constraints. We add the constraint of the maximum number of customers that can be served by a facility and denote it by $MS$. For example $MS=95$ in  $FC_{{(100,110)},{(80,90)}, MS=95}$ denotes that maximum number of customers that can be served by a facility is $95$.  In all cases, we set the number of customers and facilities for an instance to be $100$.


\subsection{Parameters}
We set number of  heads to 2 for multi-head attention with MLP hidden size 32. We train until convergence of validation loss. We set learning rate to 0.001. We use Adam optimizer for training. We set $\alpha$ and $\beta$ in eq.~\ref{eq:overall_loss}  to $1.5$ and $100$ respectively. 
\label{app:params}

\subsection{Baselines}
\label{app:baseline}

For EWC baseline(only using EWC), we set weight of the  \textit{elastic weight} component to 1000. For \GCNN we use embedding size 64.

\subsection{Number of nodes}
\begin{figure}[h!]
 \vspace{-0.20in}
\centering
\hspace{-0.5in}
\includegraphics[scale=0.38]{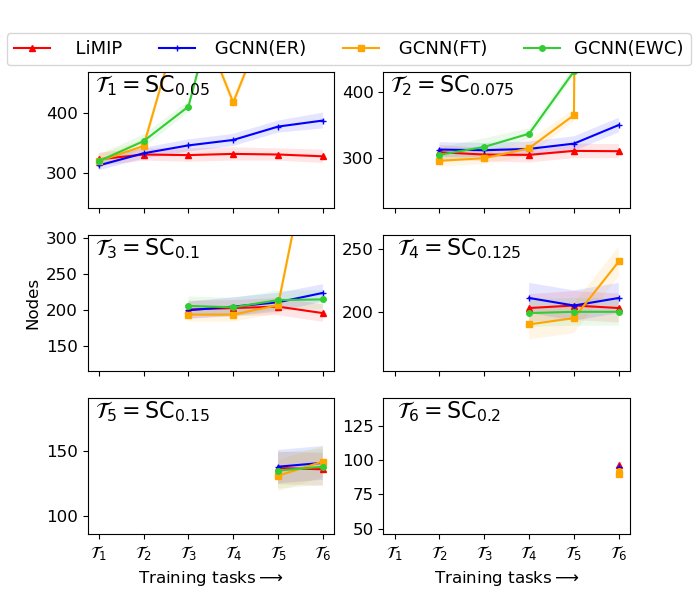}
\caption{ \textbf{Number of Nodes: Testing on SetCover in lifelong scenario}: Evolution of number of nodes  for each task when different methods are updated on each task sequentially.}
\label{fig:eval_sc_continue-node}
\end{figure}

\begin{figure}[h!]
\centering
\hspace{-0.5in}
\includegraphics[scale=0.38]{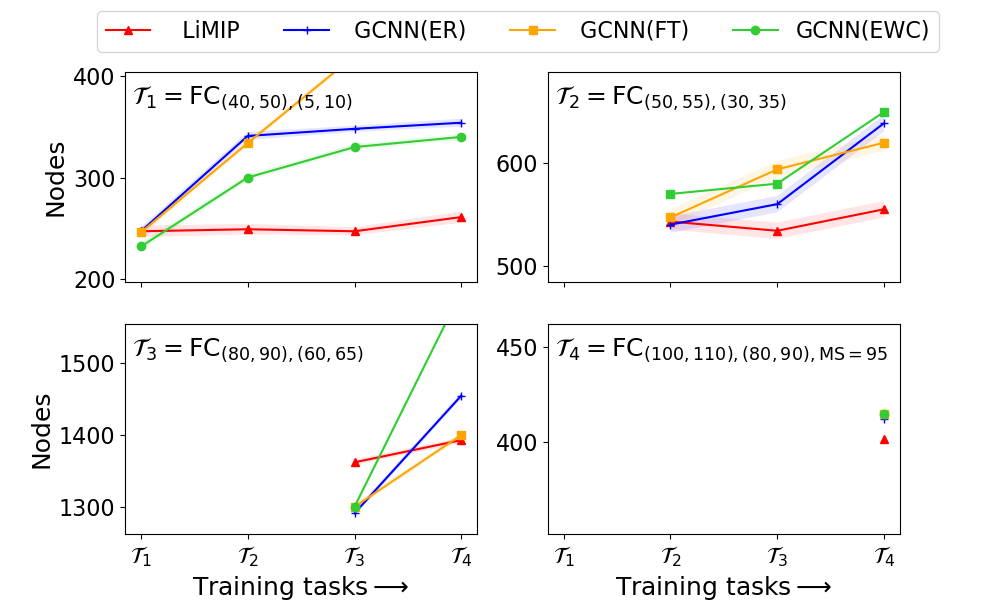}
\caption{ \textbf{Number of Nodes:  Testing on Facility location in lifelong scenario}: Evolution of number of nodes  for each task when different methods are updated on each task sequentially.}
\label{fig:eval_fac_continue-node}
\end{figure}

In addition to  plots based upon solving time shown in the main paper, in fig.~\ref{fig:eval_sc_continue-node} and ~\ref{fig:eval_fac_continue-node}, we present the number of nodes solved for different methods.
\label{app:nodes}

\subsection{Additional Results}

\label{app:results}
\subsubsection{Results on Facility location}
\begin{figure}[h!]
\centering
\hspace{-0.5in}
\includegraphics[scale=0.38]{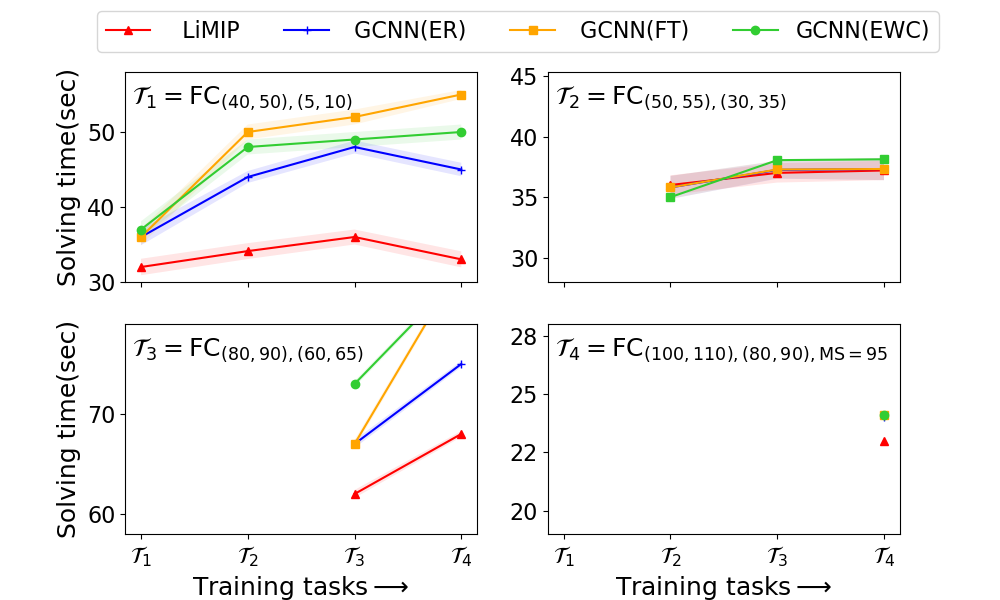}
\caption{ \textbf{Facility location: Testing on Facility location in lifelong scenario}: Evolution of solving time  for each task when different methods are updated on each task sequentially.}
\label{fig:eval_indset_continue}
\end{figure}

We present the results on Facility location problem in fig~\ref{fig:eval_indset_continue}. We observe that in facility location problem also also, \namemodel achieves significant lower forgetting in comparison to other methods. 


\subsubsection{Results on default SCIP solver heuristics}: For the sake of completeness we present results of SCIP solver in table ~\ref{tab:res:scip}. We skip results of fsb for $FC_{(40,50),(5,10)}$ due to its poor running time performance.

\begin{table}[h!]
\small
\centering
\scalebox{0.99}{
\begin{tabular}{llcc}
\toprule
\small
\textbf{Dataset} & \textbf{Method} & \textbf{Time} & \textbf{Nodes}\\ 
\midrule
\multirow{2}{*}{\textsc{$SC_{0.05}$} } & relpscost & {$9.93 \pm 1.1$} & ${497.2 \pm 16.1}$ \\
 & fsb & $63.50 \pm 2.1$  & {$91.5 \pm 4.5 $} \\
\cmidrule{1-4}

\multirow{2}{*}{\textsc{$IS_{{4},{750}}$} } & relpscost & {$66.2 \pm 0.1$} & ${21010 \pm 1.2}$ \\
 & fsb & $1531.50 \pm 0.15$  & {$1195.5 \pm0.1 $} \\
\cmidrule{1-4}

\multirow{1}{*}{\textsc{$FC_{(40,50),(5,10)}$} } & relpscost & $268.2 \pm 0.88 $  & $7528 \pm 0.62$\\
\cmidrule{1-4}

\bottomrule

\end{tabular}}%
\caption{Results of SCIP Full Strong (fsb) and SCIP Reliability pseudocost(relpscost) branching}
\label{tab:res:scip} 
\end{table}

\subsection{Poor scalability of  TreeGate~\cite{zarpellon2020parameterizing} }
\label{app:scalable}
We test the performance of TreeGate~\cite{zarpellon2020parameterizing} in our setup where we have a set of training instances available for a class of MIP. We generate relpscost branching pairs from the  \textsc{$SC_{0.05}$} dataset similar to procedure described in ~\cite{zarpellon2020parameterizing} i.e using different random seeds and also performing random branching. We used same set of training  instances as used in \namemodel and \GCNN to generate branching pairs.

We compare the running time and number of nodes.  In Table~\ref{tab:relzarpel} we observe that TreeGate is twice slower compared to imitation learning of strong branching.
\begin{table}[h!]
    \centering
        \scalebox{0.9}{
        \begin{tabular}{p{5cm}p{1.3cm}p{1.3cm}p{1.3cm}}
\toprule
        \textbf{Model} & \textbf{Time} & \textbf{\# Nodes}  \\
        \midrule
        TreeGate   &   {$24.1$} &  {$764$} \\ 
       {\namemodel} &  {$11.01$} &  {$320$}\\
       
        \bottomrule
        \end{tabular}}
    \caption{\textbf{Comparison against TreeGate\cite{zarpellon2020parameterizing}}: Comparison of Imitation learning of Strong Branching against TreeGate which learns weaker heuristic of reliability pseudo cost branching.  Both models were trained and tested on \textsc{$SC_{0.05}$}. }
    \label{tab:relzarpel}
\end{table}

\subsection{Ablation}
\label{app:ablation}

\subsubsection{\namemodel without Elastic weight consolidation}
\begin{figure}[h!]
\centering
\hspace{-0.5in}
\includegraphics[scale=0.38]{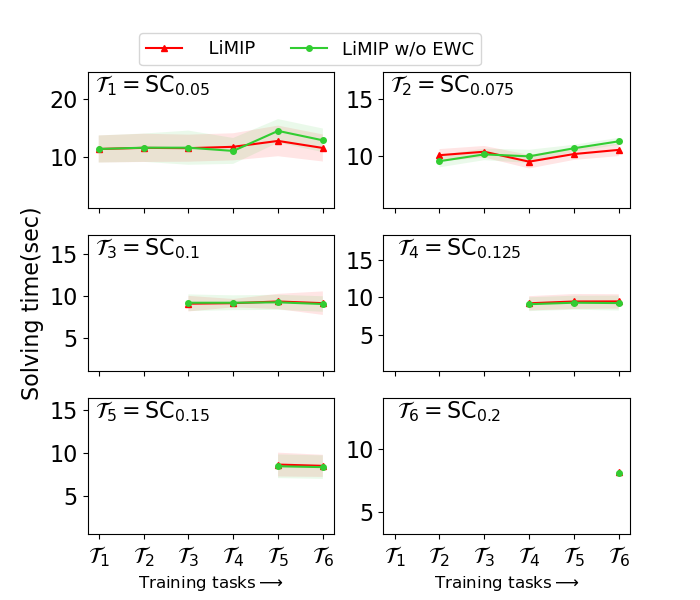}
\caption{ \textbf{\namemodel without EWC: Testing on SetCover in lifelong scenario}: Evolution of Solving time for each task when different methods are updated on each task sequentially.}
\label{fig:eval_sc_continue-time-woewc}
\end{figure}

In fig.~\ref{fig:eval_sc_continue-time-woewc} we study the performance of \namemodel without the EWC loss  $\mathcal{L}_{lifelong}$.

\subsubsection{Using small sized memory buffer:} 
\begin{figure}[h!]
\centering
\hspace{-0.5in}
\includegraphics[scale=0.38]{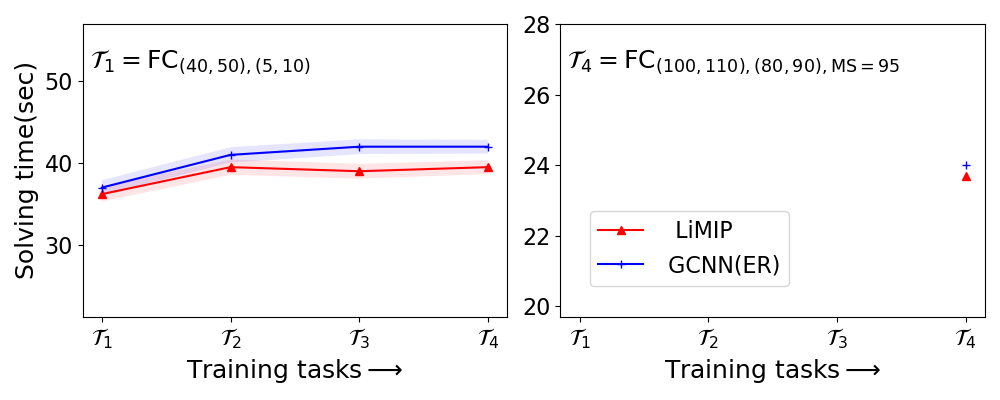}
 \vspace{-0.2in}
\caption{ \textbf{Comparison of \namemodel and ER model on  buffer size 200:} Testing on first and last problem in the Facility location problem sequence. }
\label{fig:eval_fac_buff}
\end{figure}

We use  a buffer of size $200$ instead of $500$ used earlier. We observe similar conclusion on a lower buffer size  as can be seen in fig.~\ref{fig:eval_fac_buff}. For simplicity, we show result on the first and the last problem in the facility location problem sequence.

\subsubsection{Using Bipartite GAT  encoding for baselines instead of GCNN }

\begin{figure}[h!]
\centering
\hspace{-0.5in}
\includegraphics[scale=0.38]{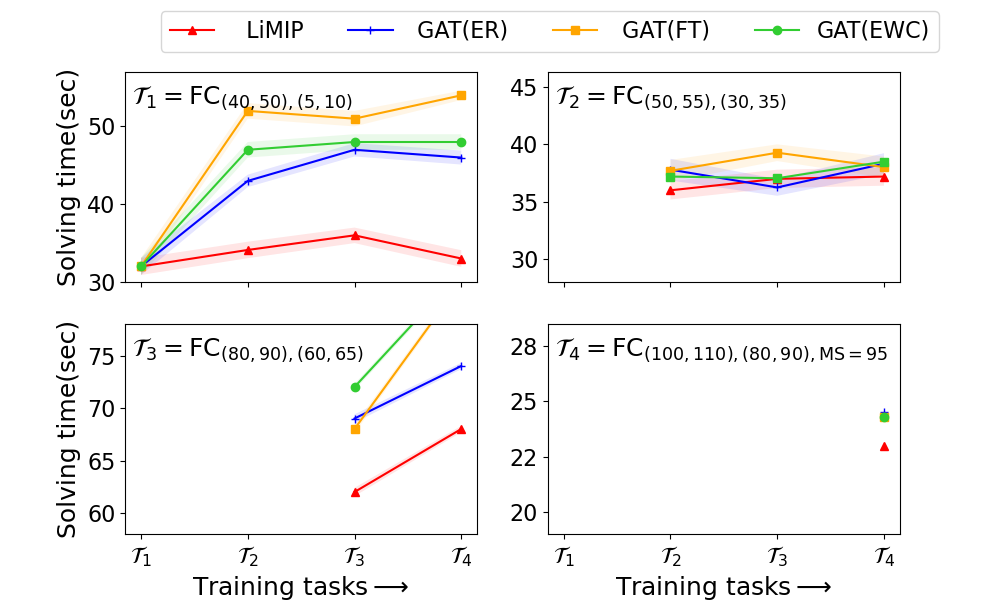}
\caption{ \textbf{Comparison using  Bipartite GAT as encoding model for baselines: Testing on Facility location in lifelong scenario}: Evolution of solving time  for each task when different methods are updated on each task sequentially.}
\label{fig:eval_fac_continue-GAT}
\end{figure}

\begin{figure}[h!]
\centering
\hspace{-0.5in}
\includegraphics[scale=0.42]{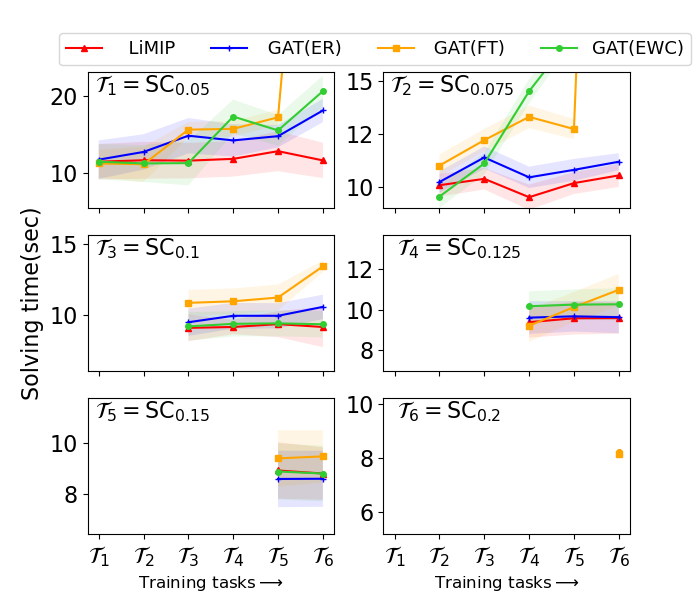}
\caption{ \textbf{Comparison using encoding model as Bipartite GAT for baselines  : Testing on SetCover in lifelong scenario}: Evolution of solving time  for each task when different methods are updated on each task sequentially.}
\label{fig:eval_sc_continue-GAT}
\end{figure}

In fig.~\ref{fig:eval_fac_continue-GAT} and ~\ref{fig:eval_sc_continue-GAT} we study the performance of different baselines when the base model used for them is bipartite GAT. This shows that \namemodel outperforms existing baselines irrespective of the GNN encoding used. We observe similar conclusion as observed on \GCNN in main paper.

%% file: 03_RelatedWork.tex
\textbf{Learning to solve MIPs}

Strong Branching is widely accepted as the most efficient branching expert in terms of number of nodes. However, its main advantage is the high computational cost of finding the variable to branch. To tackle this problem, Learning-based techniques\cite{alvarez2014supervised,alvarez2017machine,marcos2016online,khalil2016learning, gasse2019exact,nair2020solving,gupta2020hybrid} have focused on learning fast approximations of \textit{Strong Branching} rule by learning from a set of training instances for a class of MIPs. Most of these techniques either learn to rank candidate branching variables or learn to imitate expert strong branching rule. Among various learning based methods, GCNN\cite{gasse2019exact} has shown significant scalability gains and is considered the state-of-the art method for learning to branching\cite{banitalebi2021ml4co}.   It has shown to to improve upon previously proposed approaches for branching on several MILP problem benchmarks, and further, also obtains faster running time compared to the default SCIP solver. Recently, \cite{zarpellon2020parameterizing} proposed \textit{TreeGate} model to learn approximation of reliability pseudo-cost branching. However, in our study we observed, the technique does not scale well in comparison to GCNN\cite{gasse2019exact} and default SCIP solver in terms of the running time metric. Recently, there have been attempts to use reinforcement learning to obtain better heuristics\footnote{Scavuzzo, L. Learning to branch with Tree MDPs. arXiv preprint arXiv:2205.11107.}, however, in terms of running time they are still inferior to imitation learning of strong branching.

\noindent \\
\textbf{Continual Learning}

Although deep neural networks have obtained significant success on sever learning tasks, however most of them suffer from catastrophic forgetting in the continual learning scenario. The goal of continual learning is to learn to adapt to new data in a streaming scenario while consolidating the knowledge learned from previous data to prevent catastrophic forgetting. Many recent endeavours have been made towards alleviating  catastrophic forgetting. The first one's being replay-based\cite{rebuffi2017icarl,riemer2018learning,lopez2017gradient,chaudhry2018efficient} where a subset of old data samples are replayed from a memory buffer while learning a new task. Second category is the regularization based techniques\cite{zenke2017continual,kirkpatrick2017overcoming,aljundi2018memory} which learn importance of weights of the neural model for each task and prevent significant changes on them. The third category is dynamic network expansion. In this context, some progress has been made to intelligently expand the neural network\cite{yoon2017lifelong,draelos2017neurogenesis} when the current capacity of the model is not sufficient to learn new tasks without causing forgetting of earlier one's. These approaches introduced here change architectural properties in response to new information by dynamically accommodating novel neural resources with due course of time such as increased number of neurons or network layers. 
 Recently, some amount of progress has been made on continual learning on GNNs\cite{wang2020streaming,liu2021overcoming,zhou2021overcoming}. We refer the reader to the follow surveys on continual learning on neural networks\cite{parisi2019continual} and continual learning for GNNs\cite{febrinanto2022graph}.